\documentclass{amsart}
\usepackage{amssymb, amsfonts, amsmath, amsthm}
\usepackage{hyperref}
\usepackage{enumerate}
\usepackage{epsfig}
\usepackage[active]{srcltx}
\usepackage{verbatim}

\input cyracc.def

\def\RR{\mathbb{R}}%
\def\HH{\mathbb{H}}%

\def\eps{\varepsilon}
\def\theta{\vartheta}
\def\ge{\geqslant}
\def\le{\leqslant}
\def\phi{\varphi}
\def\i{\subset}

\def\<{\langle}
\def\>{\rangle}

\def\dist{\operatorname{dist}}

\def\Proof{\vspace{5mm}\noindent{\it Proof}}
\def\qeds{\qed\vspace{5mm}}

\newtheorem{thm}{Theorem}

\newtheorem{lem}[thm]{Lemma}

\newtheorem{defn}[thm]{Definition}

\newtheorem{Example}[thm]{Example}

\newtheorem{Counterexample}[thm]{Counterexample}

\newtheorem{remark}[thm]{Remark}
\newenvironment{rmk}{\begin{remark}\rm}{\end{remark}}
\newtheorem{Fact}[thm]{Fact}

\newtheorem{Nothing}[thm]{$\!\!\!$}

\begin{document}

\title[Optimal lower curvature bound]{An optimal lower curvature bound for convex hypersurfaces in Riemannian  manifolds}

\author{Stephanie Alexander}
\address{Department of Mathematics, University of Illinois at Urbana-Champaign, IL 61801} 
\email{sba@math.uiuc.edu}
\urladdr{www.math.uiuc.edu/~sba}

\author{Vitali Kapovitch}
\address{Department of Mathematics, University of Toronto,
Toronto, ON, M5S 2E4}
\email{vtk@math.toronto.edu}
\urladdr{www.math.toronto.edu/~vtk/}

\author{Anton Petrunin}
\address{Department of Mathematics, Penn State University, State College, PA 16802}
\email{petrunin@math.psu.edu}
\urladdr{www.math.psu.edu/petrunin/}


\thanks{ 
We thank the Mathematical Research Institute Oberwolfach (Research in Pairs programme)  for a valuable visit, of which this note is an offshoot. The third-named author was supported in part by the National Science Foundation under grant DMS-0406482.}

\begin{abstract}
It is proved that a convex hypersurface in a Riemannian manifold of sectional curvature $\ge \kappa$ is an Alexandrov's space of curvature $\ge \kappa$.  This theorem provides an optimal lower curvature bound for an older  theorem of Buyalo.
\end{abstract}

\maketitle

The purpose of this paper is to provide a reference for the following theorem:
\begin{thm}\label{thm:exact} Let $M$ be a Riemannian manifold with sectional curvature $\ge \kappa$. Then any convex hypersurface $F\i M$ equipped with the induced intrinsic metric is an Alexandrov's space with curvature $\ge \kappa$.
\end{thm}

Here is a slightly weaker statement:

\begin{thm} \cite{buyalo:convex-surface}\label{buy}
If $M$ is a Riemannian manifold, then any convex hypersurface $F\i M$ equipped with the induced intrinsic metric is locally an Alexandrov's space.
\end{thm}

In the proof of Theorem \ref{buy} in \cite{buyalo:convex-surface}, the (local) lower curvature bound depends on (local) upper as well as lower curvature bounds of $M$.
We show that the approach in \cite{buyalo:convex-surface}
can be modified  to give Theorem \ref{thm:exact}.

\begin{defn} \label{def:lambda-conv}A locally Lipschitz function $f$ on an open subset of a Riemannian manifold is called \emph{$\lambda$-concave} ($\lambda \in \RR$) if for any unit-speed geodesic $\gamma$, 
the function 
$$f\circ\gamma(t) - \lambda t^2/2$$
is concave.
\end{defn}

\begin{lem}\label{lem:smoothing}
Let  $f:\Omega\to\RR$ be a $\lambda$-concave function on an open subset $\Omega$ of a Riemannian manifold. 
Then there is a sequence of nested open domains $\Omega_i$, with $\Omega_i\i\Omega_j$ for $i<j$ and $\cup_i\Omega_i=\Omega$, and a sequence of smooth $\lambda_i$-concave functions $f_i:\Omega_i\to\RR$ such that
\begin{enumerate}[(i)]
\item on any compact subset $K\i\Omega$, $f_i$ converges uniformly to $f$;
\item $\lambda_i\to\lambda$ as $i\to\infty$.
\end{enumerate}
\end{lem} 

This lemma is a slight generalization of \cite[Theorem 2]{greene-wu} and can be proved exactly the same way.

\Proof\ \textit{of Theorem \ref{thm:exact}}. 
Without loss of generality one can assume that 
\begin{enumerate}[(a)]
\item \label{k>=-1}$\kappa\ge -1$, 
\item $F$ bounds a compact convex set $C$ in $M$, 
\item there is a $(-2)$-concave function $\mu$ defined in a neighborhood of $C$ and $|\mu(x)|<1/10$ for any $x\in C$,
\item \label{unique} there is unique minimal geodesic between any two points in $C$. 
\end{enumerate}
(If not, rescale and pass to the boundary of the convex piece cut by $F$  from a small convex ball centered at $x\in F$, taking $\mu=-10\dist_x^2$.)

Consider the function $f=\dist_F$.
By Rauch comparison (as in \cite[11.4.8]{petersen}), for any unit-speed geodesic $\gamma$ in the interior of $C$, $(f\circ\gamma)''$ is bounded in the support sense by the corresponding value in the model case (where $M=\HH^2$ and $F$ is a geodesic).  In particular,
$$(f\circ\gamma)''\le f\circ\gamma.$$

Therefore $f+\eps \mu$ is $(-\eps)$-concave in 
$\Omega_\eps=C\cap f^{-1}((0,\eps))$.
Take $K_\eps=f^{-1}([\frac13\eps,\frac23\eps])\cap C$. 
Applying lemma \ref{lem:smoothing}, we can find a smooth $(-\tfrac\eps2)$-concave function $f_{\eps}$ which is arbitrarily close to $f+\eps \mu$ on $K_\eps$ and which is defined on a neighborhood of $K_\eps$. 
Take a regular value $\theta_\eps\approx\tfrac12\eps$ of $f_\eps$. (In fact one can take $\theta_\eps=\tfrac12\eps$, but it requires a little work.) 
Since $|\mu|_C|<1/10$, the level set $F_\eps=f_{\eps}^{-1}(\theta_\eps)$ will lie entirely in $K_\eps$.
Therefore $F_\eps$ forms a smooth closed convex hypersurface.
By the Gauss formula, the sectional curvature of the induced intrinsic metric of $F_\eps$ is $\ge\kappa$.
$F_\eps$ bounds a compact convex set $C_\eps$, where $F_\eps\to F$, $C_\eps\to C$ in Hausdorff sense as $\eps\to 0$. By property (\ref{unique}), the restricted metrics from $M$ to $C, C_\eps$ are intrinsic, and so $C_\eps$ is an Alexandrov space with $F_\eps$ as boundary, that converges in Gromov--Hausdorff sense to $C$.  It follows from \cite[Theorem 1.2]{petrunin:extremal} (compare \cite[Theorem 1]{buyalo:convex-surface}) that $F_\eps$ equipped with its intrinsic metric converges in Gromov--Hausdorff sense to $F$ equipped with its intrinsic metric.
Therefore $F$ is an Alexandrov space with curvature $\ge \kappa$.\qeds

\begin{rmk}
We are not aware of any proof of theorem \ref{thm:exact} which is not based on the Gauss formula. (Although if $M$ is Euclidean space, there is a beautiful purely synthetic proof in \cite{milka-conv}.) 
Finding such a proof would be interesting on its own, and also could lead to the generalization of theorem \ref{thm:exact} to the case when $M$ is an Alexandrov space.
\end{rmk}

\end{document}